\documentclass[12pt]{amsart}
\usepackage{latexsym}
\relax
\catcode`"\active
\bibstyle{plain}
\begin{document}

\title{
LEIBNIZ'S DEFINITION OF MONAD
}
\author{S. S. Kutateladze}
\address[]
{Sobolev Institute of Mathematics\newline
\indent Siberian Division of the\newline
\indent Russian Academy of Sciences\newline
\indent Novosibirsk, RUSSIA
}
\email{
sskut@member.ams.org
}
\begin{abstract}
This is a short discussion of the  definition of monad
which was given by G.~W. Leibniz in his {\it Monadology}.
\end{abstract}
\date{August 11, 2006}

\maketitle

Acquired traits are never inherited. This law of genetics
determines many aspects of public life. Mankind
creates and supports complicated social institutions
for transferring to the young generations the experience of their
ancestors.  As biological species,
we differ little from our paleolithic predecessors.
So we may hope to comprehend the thoughts and ideas
that are bequeathed to us by the greatest minds of the past epochs.

The outlook of Leibniz, proliferating with his works,
occupies a unique place in  human culture.
We can hardly find in the philosophical treatises of his
predecessors and later thinkers something
comparable with the phantasmagoric conceptions of monads,
the special and stunning
constructs of the world and mind which precede, comprise, and
incorporate all the infinite advents of the eternity.
{\it Monadology} \cite{Mon} is usually dated as of 1714.
This article was never published during Leibniz's life.
Moreover, it is generally accepted  that the very term ``monad''
had appeared in his writings since 1690 when he was already an
established and prominent scholar.

The special attention to the origin of the
term ``monad'' and the particular investigation into
the date of its first appearance in the works by Leibniz are in
fact the present-day products.
There are now a few if any cultivated persons who never
got acquaintance with the basics of planimetry and heard nothing of
Euclid. However, no one has ever met the concept of ``monad''
on the school bench.
Neither the contemporary translations of Euclid's {\it Elements}
nor the popular school text-books
contain this seemingly exotic term. However,
the concept of  ``monad''  is fundamental not only
for Euclidean geometry but also for the whole science of
the Ancient Hellada.

By Definition I of Book VII of Euclid's {\it Elements} \cite{Euc}
a~monad is ``that by virtue of which each
of the things that exist is called one.''
Euclid proceeds with Definition~2:
``A number is a multitude composed of monads.''
Note that the present-day translations of the Euclid treatise
substitute ``unit'' for ``monad.''

A contemporary reader can hardly understand why
Sextus Empiricus, an outstanding scepticist
of the second century, wrote when
presenting the mathematical views of his predecessors
as follows \cite{Sextus}:
``Pythagoras said that the origin of the
things that exist is a~monad by virtue of which
each of the things that exist is called one.''
And furthermore:
``A~point is structured as a~monad; indeed,  a~monad is
a~certain origin of numbers and likewise a~point is a~certain origin
of lines.''
Now some place is in order for the excerpt which
can easily be misconceived as a citation
from {\it Monadology}:
``A~whole as such is indivisible and
a~monad, since it is a~monad, is not divisible.
Or, if it splits into many pieces
it becomes a~union of many monads
rather than a~[simple] monad.''

It is worth observing that the ancients
sharply perceived an exceptional status of the
start of counting. In order to count, one
should firstly  particularize the entities to count and
only then to proceed with putting these entities into correspondence
with some symbolic series of numerals.
We begin counting with making ``each of the things one.''
The especial role of the start of counting
is reflected in the almost millennium-long dispute
about whether or not the unit (read, monad) is a natural number.
We feel today that it is excessive to distinguish the key role
of the unit or monad which signifies the start of counting.
However, this was not always so.

From the times of Euclid, all serious scientists knew
about existence of the two basic concepts of mathematics: a point and a monad.
By Definition~1 of Book~1 of Euclid's {\it Elements}:
``A point is that which has no parts.''
Clearly this definition differs drastically from the definition of
monad as that which makes one from many.
The cornerstone  of geometry is other than that of arithmetic.
Without clear understanding of this circumstance it is impossible
to comprehend the essence of the views of Leibniz.
By the way, the modern set theory refers to ``that which has no parts''
as the empty set, the starting cardinal of the von Neumann universe.
The present-day mathematics seems to have no concept
that is vocalized as ``that which many makes into one.''
We will return to the modern mathematical definition of monad
shortly.

Attempting to pursue the way of Leibniz's thought,
we must always keep in mind that he was a mathematician by
belief. From his earliest childhood,  Leibniz dreamed of
``some sort of calculus'' that operates in the
``alphabet of human thoughts'' and possesses
the same beauty, strength, and integrity as mathematics in solving
arithmetical and geometrical problems. Leibniz devoted many articles
to invention of this universal logical calculus.
The diversity and even polarity of the views of these
writings proceed along with the universally accepted appraisal of Leibniz
as a key figure of the prehistory of the modern mathematical logic.
{\it Monadology} is listed alongside the classical achievements of Leibniz which we
express with the words {\it culculamus} and {\it differentia}.

Leibniz always emphasized his love  and devotion to mathematics.
He  stressed constantly that his general methodological views
base on ``study into the methods of analysis in mathematics
which I was engrossed in with such an eager that I do not know
whether it is possible to find many who served it with more toil.

As a top mathematician of his age, Leibniz
was in full command of Euclidean geometry.
Therefore, we are upmost bewildered  already to
read Item~1 of his {\it Monadolody} where he gave the first
impression about  his monad:
``The Monad, of which we shall here speak, is nothing but a
simple substance, which enters into compounds. By 'simple' is meant
'without parts.'''
This definition of monad as a ``simple'' substance
without parts coincides with the Euclidean definition of
point. At the same time the reference to compounds consisting of monads
reminds us the structure of the definition of number which belongs to Euclid.

The synthesis of both primary definitions of Euclid
in the Leibnizian monad is not accidental.
We must always bear in mind that the seventeenth century
is the epoch of microscope. It was already in the 1610s that microscopes
were mass-produced in many European countries. From the 1660s
Europe was enchanted by Antony van Leeuwenhoek's microscope.

Let us make a mental experiment
and aim  a strong microscope at
a region about  a point at a mathematical line.
We will see in the
eyepiece a~blurred and dispersed cloud with unclear frontiers
which is a visualization of the point under investigation.
Under greater magnification,
the portion of the ``point-monad'' we are looking at
will enlarge, revealing extra details whereas
disappearing  partially from sight.
However, we are still inspecting the same standard real number
which  you might prefer to percept as
described by this process of
``studying the microstructure of a~physical straight line.''
Visualizing a point by microscope reveals
its monadic essence. Leibniz could reason so or approximately so.
In any case, the view of the monad of a standard real number
as the collection of all infinitely close points
is generally adopted in the contemporary infinitesimal
analysis resurrected under the name of {\it nonstandard analysis}
in the works by Abraham Robinson in 1961.

\bibliographystyle{plain}

\end{document}